\documentclass[12pt]{article}

\usepackage{graphicx}
\usepackage{float}

\begin{document}

\title{On the statistics of differences of zeta zeros starting from zero number $10^{23}$}

\maketitle
\begin{center}
\author{Jouni Takalo}

{Dpt. of Physical Sciences, Space Climate research unit, University of Oulu,
POB 3000, FIN-90014, Oulu, Finland\\
email:{jouni.j.takalo@oulu.fi}}

\end{center}

\begin{abstract}
We study distributions of differences of unscaled Riemann zeta zeros, $\gamma-\gamma^{'}$, at large. We show, that independently of the location of the zeros, i.e., even for zeros as high as $10^{23}$, their differences have similar statistical properties. The distributions of differences are skewed usually towards the nearest zeta zero. We show, however, that this is not always the case, but depends upon the distance and number of nearby zeros on each side of the corresponding distribution. The skewness, however, always decreases when zeta zero is crossed from left to right, i.e., in increasing direction. Furthermore, we show that the variance of distributions has local maximum or, at least, a turning point at every zeta zero, i.e., local minimum of the second derivative of the variance. In addition, it seems that the higher the zeros the more compactly the distributions of the differences are located in the skewness-kurtosis -plane. Furthermore, we show that distributions can be fitted with Johnson probability density function, despite the value of skewness or kurtosis of the distribution.
\end{abstract}

 \textbf{Keywords:} Riemann nontrivial zeta zeros, difference of zeta zeros, distribution function, skewness, variance, Johnson probability distribution

\section{Introduction}

Montgomery (1973) conjectured that the paircorrelation between pairs of zeros of the Riemann zeta function (scaled to have unit average spacing) is \cite{Montgomery}

\begin{equation}
\label{eq:CUE}
	R_{2}\left(x\right) = 1\!-\!\left(\frac{sin(\pi\,x)}{\pi\,x}\right)^{2} .
\end{equation}

Odlyzko (1987) already showed that this conjecture was supported also by larger heights of zeta zeros \cite{Odlyzko}. The aforementioned papers are restricted to analyse only consecutive or locally close differences of zeros. We study here distributions of differences of unscaled zeta zeros at large. Perez-Marco (2011) has shown that the statistics of very large zeros do find the location the first Riemann zeros \cite{Perez-Marco}. Takalo (2020) showed that distributions of differences of zeros are skewed towards the nearest zeta zero, and have local maximum of variance at or near zeta zero \cite{Takalo}. Here we study distributions of differences starting from $10^{23}$-rd zero, and show that the properties even at that height have similar properties and are aware of the very first zeta zeros, i.e., information about low zeros is encoded also to large zeta zeros and their differences. Furthermore, we use Johnson probability density function (PDF) to show its flexibility for fitting zeta zero difference distributions.

\section{Data and methods}

\subsection{Nontrivial zeros of Riemann zeta function}

The data, i.e., the imaginary parts of nontrivial Riemann zeta zeros were fetched from  (https://www.lmfdb.org/zeros/zeta/) for the five million zero starting from 1 billionth and 100 billionth zero. The zeros starting from zero number $10^{23}$ were kindly provided by Dr. A. Odlyzko.

\subsection{Johnson distribution}

Johnson distribution for the variable x is defined as 
\begin{equation}
z = \lambda+\delta\,ln\left(f\left(u\right)\right),
\end{equation}
with 
\begin{equation}
u = \left(x-\xi\right)/\lambda,
\end{equation}
Here z is a standardized normal variable and $f\left(u\right)$ has three different forms
the lognormal distribution, $S_{L}$:
\begin{equation}
f\left(u\right)=u,
\end{equation}
the unbounded distribution, $S_{U}$:
\begin{equation}
f\left(u\right)=u+{\left(1+u^{2}\right)}^{1/2},
\end{equation}
and the bounded distribution, $S_{B}$:
\begin{equation}
f\left(u\right)=u/\left(1-u\right).
\end{equation}
The supports for the distributions are $S_{L}: \xi<x,\; S_{U}: - \infty<x< \infty\; $and $S_{B}: \xi<x<\xi+\lambda$ \cite{Johnson, Wheeler}. The reason for the transformation of the the non-normal variables to standardized normal variables was that normal distribution was the only well-defined distribution at those times. However, with these definitions, the probability distributions are for
$S_{L}$:
\begin{equation}
P\left(u\right)=\frac{\delta}{\sqrt{2\pi}}\,\times\,\frac{1}{u}\,\times\,exp\left\{-\frac{1}{2}\left[\gamma+\delta\ln\left(u\right)\right]^{2}\right\}.
\end{equation}
for $S_{U}$
\begin{equation}
P\left(u\right)=\frac{\delta}{\sqrt{2\pi}}\,\times\,\frac{1}{\sqrt{u^{2}+1}}\,\times\,exp\left\{-\frac{1}{2}\left[\gamma+\delta\ln\left(u+\sqrt{u^{2}+1}\right)\right]^{2}\right\}.
\end{equation}
and for $S_{B}$
\begin{equation}
P\left(u\right)=\frac{\delta}{\sqrt{2\pi}}\,\times\,\frac{1}{u/\left(1-u\right)}\,\times\,exp\left\{-\frac{1}{2}\left[\gamma+\delta\ln\left(\frac{u}{1-u}\right)\right]^{2}\right\}
\end{equation}

\section{Distributions of differences of zeta zeros}

In the next we study the differences, $(delta)$, of unscaled zeta zeros. We use the following notation

\begin{equation}
\label{eq:delta}
\delta(n) = \gamma(n+i)-\gamma(i), n=1,2,3,...
\end{equation}

\begin{figure}
	\centering
	\includegraphics[width=0.95\textwidth]{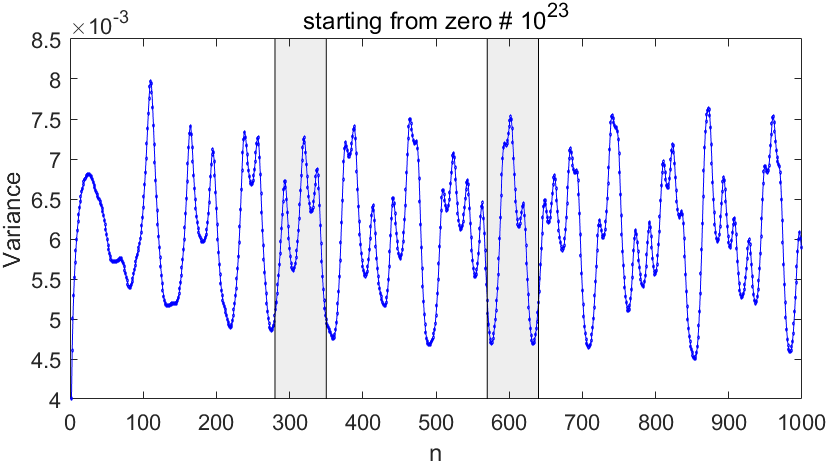}
		\caption{Variances of the $delta(n)$-distributions with n=1,2,...,999 calculated with 5 million unscaled zeta zeros starting from zero $\# 10^{23}$.}
		\label{fig:Variances_n}
\end{figure}

Here $i$ goes from j to j+5000000 in our analyses, where j is the ordinal number of the starting zeta zero. We study mainly the differences of the zeros at height $10^{23}$, and only compare those to the statistical properties of lower level (1 billion and 100 billion) differences of zeta zero. Figure \ref{fig:Variances_n} shows the variances of 5 million differences starting from zero $\# 10^{23}$ as a function of n in Eq.\ref{eq:delta} \cite{Takalo}. 

\begin{figure}
	\centering
	\includegraphics[width=1.0\textwidth]{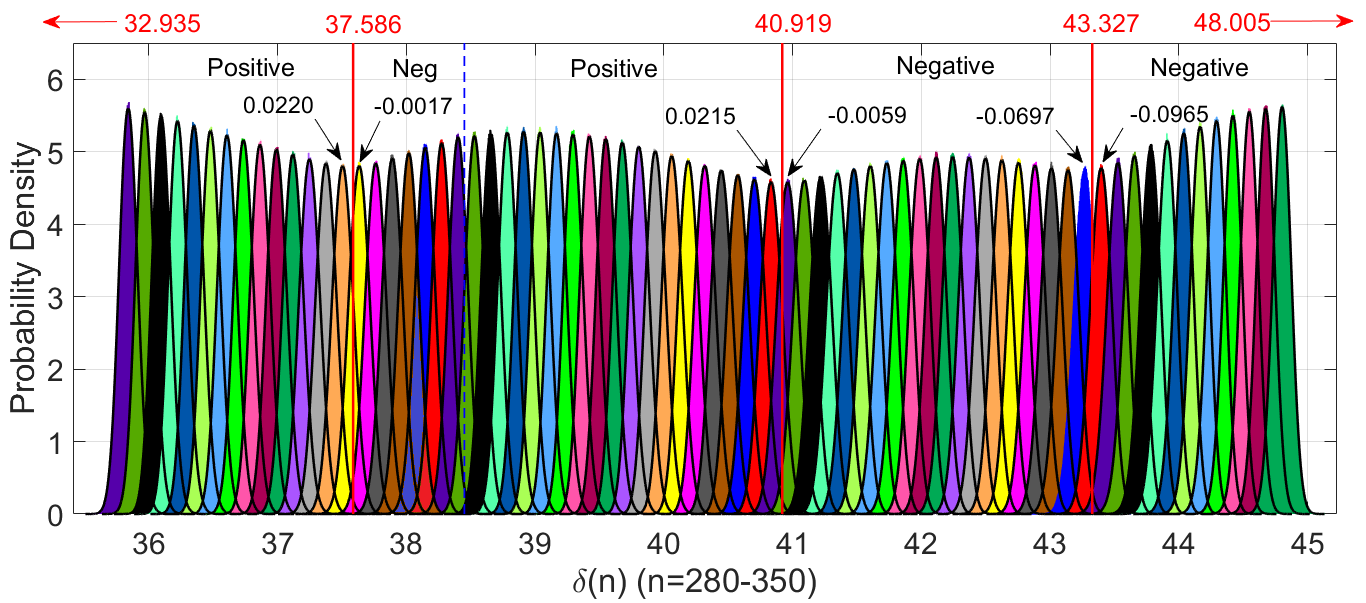}
		\caption{$Delta(n)$-distributions calculated with 5 million unscaled zeta zeros with n=280-350. Distributions are fitted with Johnson distributions.}
		\label{fig:Delta_n_280_350}
\end{figure}

Figure \ref{fig:Delta_n_280_350} shows the distributions of differences calculated with n=280-350 of the first strip marked as light gray in Fig. \ref{fig:Variances_n}. The distributions are fitted with the Johnson PDF. The skewnesses (shown for some distributions with arrowed number) of the distributions start with positive (right-hand) skewness until the zeta zero at 37.586. After the zero skewness changes to negative (left-hand) skewness, i.e., again towards the zeta zero. The skewness changes again to positive at the site marked with dashed blue vertical line. Note, that the change is not symmetrically in the middle of the zeros 37.586 and next zero at 40.919. This is because there are two nearby zeros on the right side and only one nearby zero in the left side (the second left side zero is quite far away, i.e., at 32.935 shown with red arrow). The skewness changes again at the zero 40.919 to negative and stays negative to the end of the considered interval. One can ask why the skewness does not change sign at zero 43.327. This is probably, because the next zero on the right side is farther away than the zeros in the left side. Anyway, the skewness of distribution always decreases when passing from left to right side of the zero.

\begin{figure}
	\centering
	\includegraphics[width=1.0\textwidth]{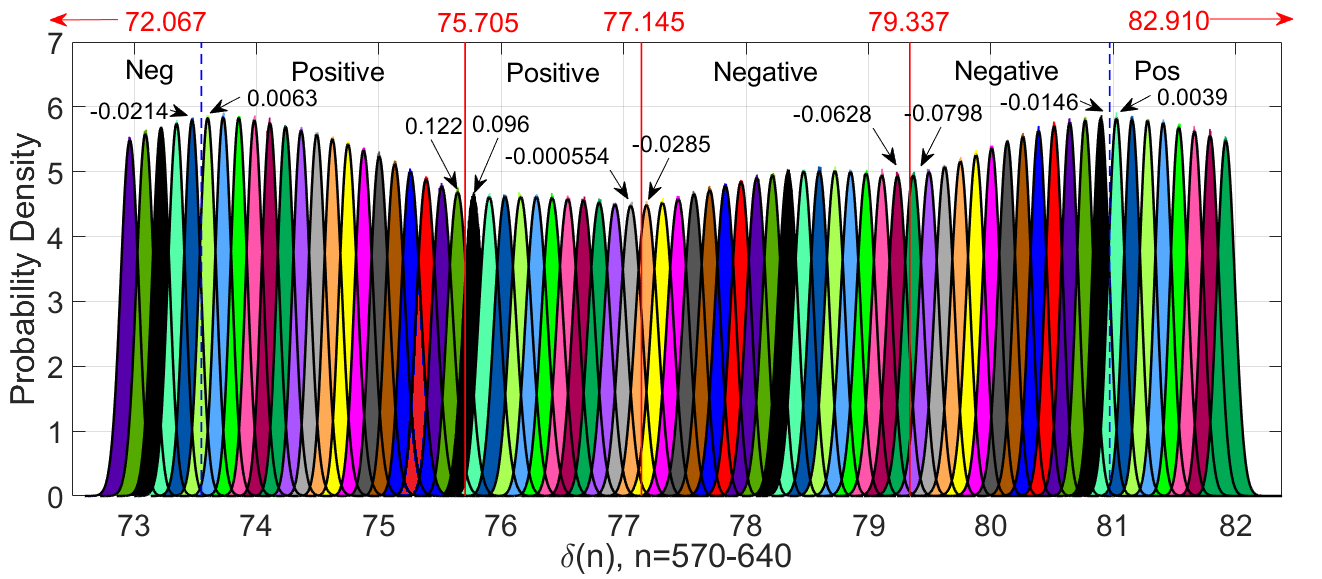}
		\caption{$Delta(n)$-distributions calculated for 5 million unscaled zeta zeros with n=570-640.  Distributions are fitted with Johnson distributions.}
		\label{fig:Delta_n_570_640}
\end{figure}

Figure \ref{fig:Delta_n_570_640} shows another interval of the distributions for n=570-640 (second strip in Fig. \ref{fig:Variances_n}). The skewnesses of the distributions start with negative skewness because of the nearby zeta zero at 72.067 (not seen in the figure). At blue dashed line skewness changes from negative to positive, because of the next zero 75.705. Interestingly, after the zero skewness does not change to negative only decreases somewhat. This is because the two next zeros at right side 77.145 and 79.337 are together more powerful than the only zero on the left side. The skewness changes sign then at 77.145, except that the last distribution at the left side is already slightly negative. Note that skewness stays negative when crossing the zero at 79.337, only decreasing again somewhat, until the blue dashed line, after which the lurking zero at 82.910 (not seen in the figure) turns the distributions to positively skewed again.

\begin{figure}
	\centering
	\includegraphics[width=1\textwidth]{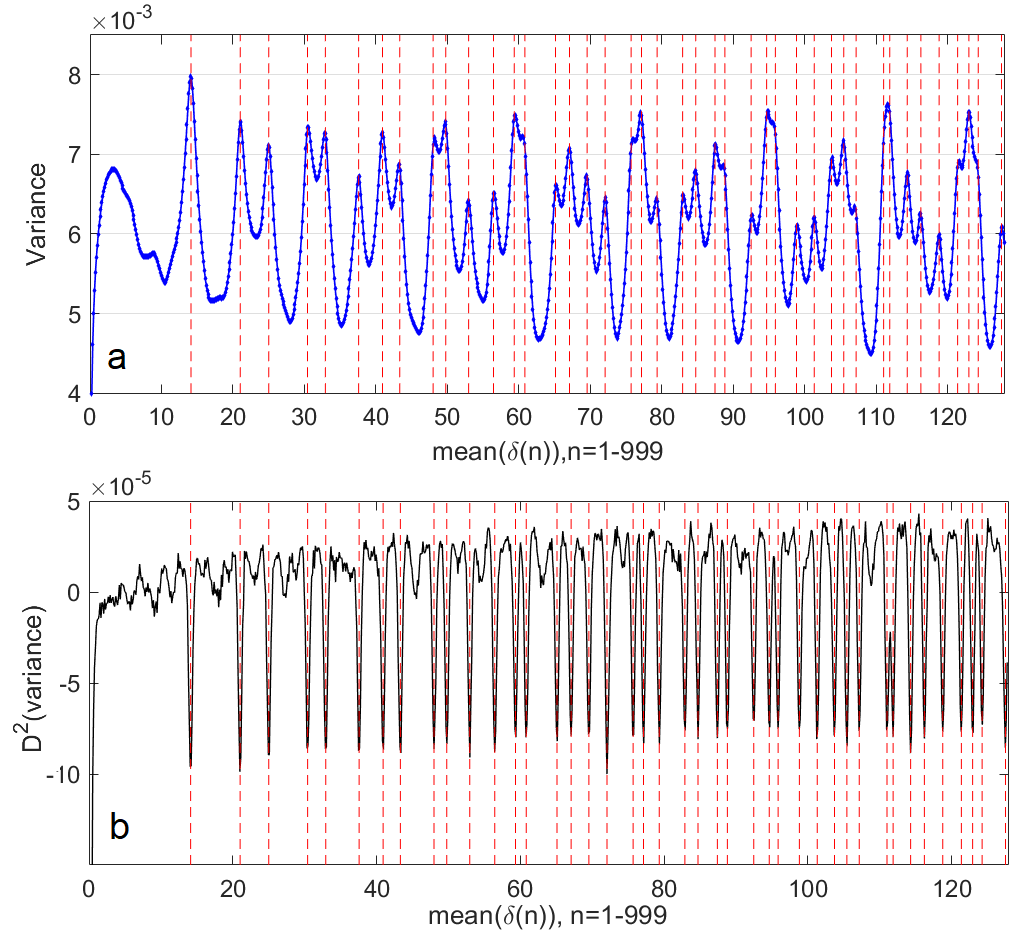}
		\caption{a) Variances of the $delta(n)$-distributions with n=1,2,...,999 calculated with 5 million unscaled zeta zeros starting from zero $\# 10^{23}$ as a function of mean($delta(n)$). b) Second derivative of the variance-curve of figure a. Zeta zeros are shown with vertical dashed red lines.}
		\label{fig:Variances_n_1_999}
\end{figure}

\begin{figure}
	\centering
	\includegraphics[width=1\textwidth]{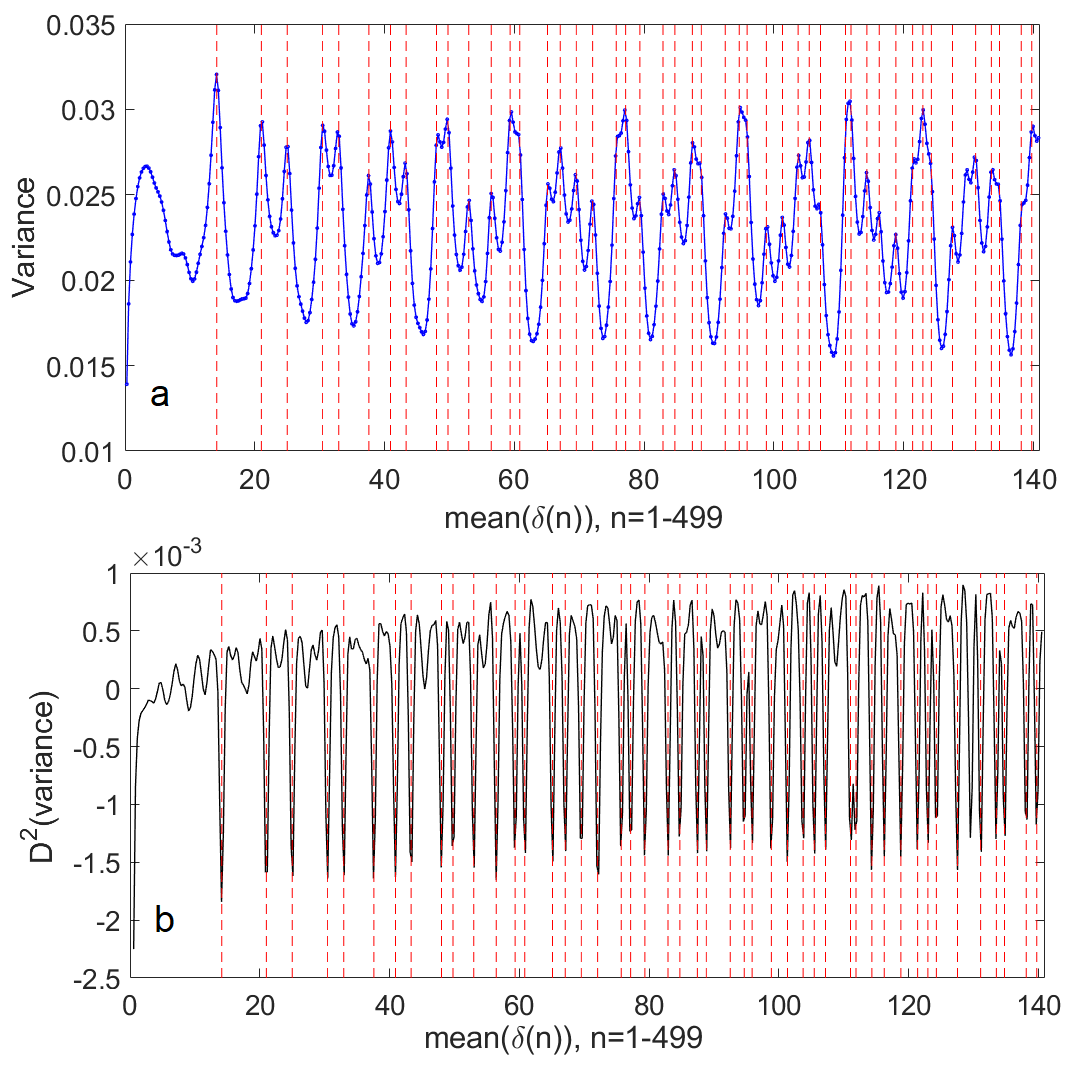}
		\caption{a) Variances of the $delta(n)$-distributions with n=1,2,...,499 calculated with 5 million unscaled zeta zeros starting from zero 100 billion as a function of mean($delta(n)$). b) Second derivative of the variance-curve of figure a. Zeta zeros are shown with vertical dashed red lines.}
		\label{fig:Variances_100_billion} 
\end{figure}

Figure \ref{fig:Variances_n_1_999}a shows variances of the distributions with n=1,2,..,999 (each separate variance shown with dots in the curve) and the Fig. \ref{fig:Variances_n_1_999}b second derivative of the variances of Fig. \ref{fig:Variances_n_1_999}a as a function of mean values of the distributions $delta(n)$. Note, that variance has local maximum or, at least, a turning point at every zeta zero. This is more clearly seen in Fig. \ref{fig:Variances_n_1_999}b as a minimum of second derivative at each zeta zero. It, indeed, seems that the first zeros are encoded to the zeros (or their differences) even at height $10^{23}$.
For comparison, we show variances and second derivative of variances for differences of zeta zeros at 100 billionth zero in Figs. \ref{fig:Variances_100_billion}a and b, respectively. The Figs. \ref{fig:Variances_n_1_999} and \ref{fig:Variances_100_billion} are very similar, although the variances of the former figure are calculated from differences of zeros twice as high than the variances of the latter figure. The difference is in the scales, which is smaller in vertical direction of Fig. \ref{fig:Variances_n_1_999}. The distribution are also narrower at height $10^{23}$ and consequently more than twice points of variances are needed to reach the same information about the zeros at height $10^{23}$ than at height 100 billion. Note also, that the curve of the second derivative of variances of Fig. \ref{fig:Variances_n_1_999}b is somewhat more noisy between the zero lines than the curve of the second derivative of Fig. \ref{fig:Variances_100_billion}b.

\begin{figure}
	\centering
	\includegraphics[width=1\textwidth]{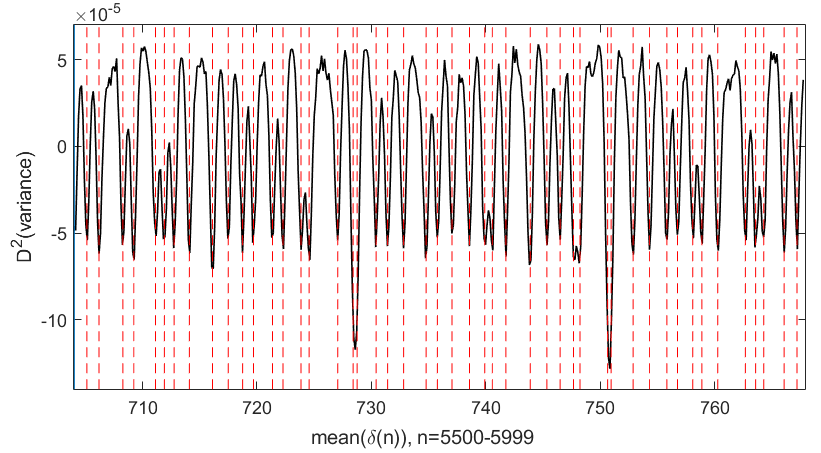}
		\caption{Second derivative of the variances of the $delta(n)$-distributions with n=5500-5999 calculated for 5 million unscaled zeta zeros starting from zero $\# 10^{23}$ as a function of mean($delta(n)$). Zeta zeros are shown with vertical dashed red lines.}
		\label{fig:Variances_n_5500}
\end{figure}

Figure \ref{fig:Variances_n_5500} shows the second derivative of the variances in the interval 704-767, i.e., for n=5500-5999 with the 48 zeta zeros belonging to the same interval. The information about these somewhat upper zeros is also included in the differences at height $10^{23}$, although resolution is too poor to distinguish the two, very nearby, pairs (728.405, 728.759) and (750.656, 750.966) as separate peaks.

\begin{figure}
	\centering
	\includegraphics[width=1.0\textwidth]{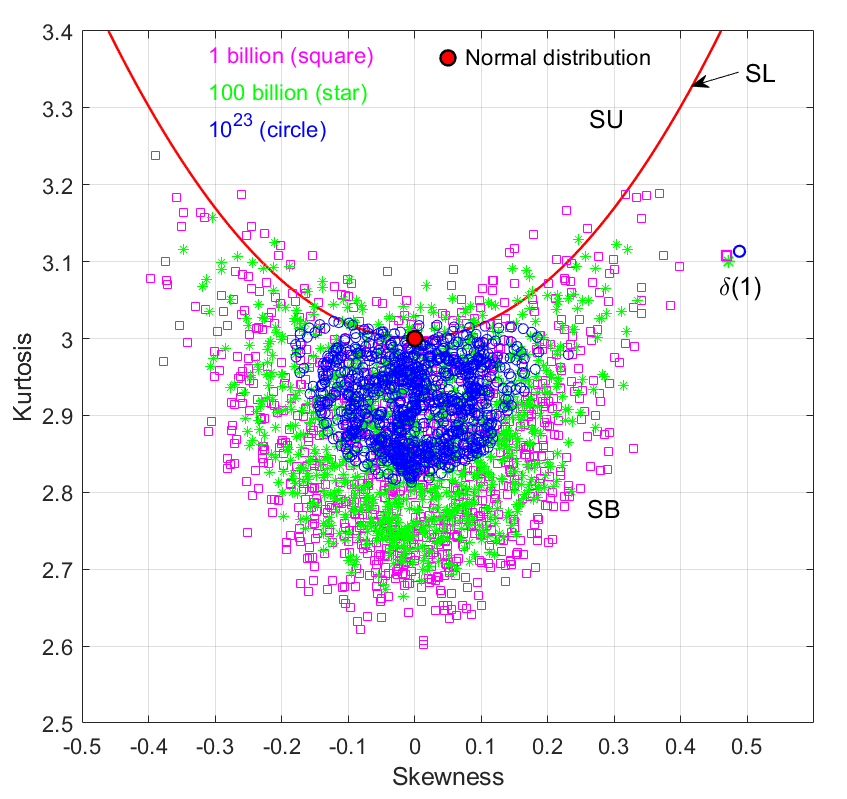}
		\caption{Skewness-kurtosis -plane with points of distributions starting 1 billionth (magenta square), 100 billionth (green star) and $10^{23}$ -rd (blue circle) zero. The red curve is Johnson SL distribution, which divides the plane to Johnson SU and SB regions. Red dot at (0,3) is the site of normal distribution in the skewness-kurtosis -plane. }
		\label{fig:Skewness_kurtosis_coord} 
\end{figure}

As shown earlier the $delta$-distributions can be fitted very well with the Johnson SB and SU PDFs (very rarely SL is needed). Figure \ref{fig:Skewness_kurtosis_coord} shows the distributions in the skewness-kurtosis -plane \cite{Cugerone} for zeros at 1 billion, 100 billion and $10^{23}$. The border between Johnson SU and SB distributions is marked with red curve in the plane (this is also region of SL distribution). It is notable that $delta(1)$s are almost in the same place in the plane for all groups. Otherwise, the points seem to be more compactly located when going to higher levels of zeta zeros. Note, that all groups of distributions with n=2-999 form a heart-shaped pattern mainly below the point (0,3), which is the site of normal distribution in the skewness-kurtosis -plane. We suppose that the average kurtosis is approaching the value 3, when going still to upper levels of zeta zeros. On the contrary, the average skewness may stay somewhat positive, because there are always new zeros lurking in the right side of the distributions, at least, if the Riemann-hypothesis is true.

\newpage

\section{Conclusions}

We have studied $delta$-distributions of zeros of Riemann zeta function at heights 100 billionth and $10^{23}$ -rd zero such that we calculate distribution for each difference, $delta(n)$, separately. We used 5 million $delta$s for these analyses, and showed that statistical properties are very similar for both intervals. The skewness of the distributions changes sign, when crossing zeta zero or, at least, decreases when moving from left side to right side of the zero (see also \cite{Takalo}). In addition, the variance has local maximum or turning-point at zeta zero. In both cases the second derivative of the variance has local minimum at each zeta zero.
We also plotted $delta(n)$-distributions (for n=1-999) of zeta zeros at heights 1 billion, 100 billion and $10^{23}$ on the skewness-kurtosis -plane. All these groups of points form a heart-shaped pattern such that the higher the zeros the more compactly the corresponding points are located in the skewness-kurtosis -plane. We believe, that going still higher levels in the zeros the average kurtosis approaches value 3, which is the kurtosis of normal distribution. $Delta$(1) is located almost at the same point for all groups apart from other points of the groups. 
\newline

\textbf{Acknowledgements:}

We acknowledge LMFDB and A.M. Odlyzko for the zeta zero data.

\end{document}